\documentclass[nonatbib]{elsarticle}

\makeatletter
\let\c@author\relax
\makeatother

\usepackage[utf8]{inputenc}
\usepackage[style=apa, citestyle=authoryear-comp, uniquelist=false, sorting=nyt, backend=biber, uniquename=false, hyperref=true, sortcites=false]{biblatex}
\usepackage{hyperref}
\usepackage{mathptmx}
\usepackage{mathtools}
\usepackage{amsmath} 
\usepackage{braket}
\usepackage{float}
\usepackage{lipsum}  
\usepackage{svg}

\usepackage{mathptmx}
\usepackage{stackrel}
\usepackage{mathtools}
\usepackage{amsmath} 
\usepackage{braket}
\usepackage{float}
\usepackage{lipsum}  
\usepackage{svg}
\usepackage{soul}

\begin{document}

\title{Queues with service resetting}

\author[1]{Ofek Lauber Bonomo}
\ead{ofekzvil@mail.tau.ac.il}

\author[2]{Uri Yechiali}
\ead{uriy@tauex.tau.ac.il}

\author[1]{Shlomi Reuveni\corref{cor1}}
\ead{shlomire@tauex.tau.ac.il}

\cortext[cor1]{Corresponding author}

\affiliation[1]{organization={School of Chemistry, Center for the Physics \& Chemistry of Living Systems, Ratner Institute for Single Molecule Chemistry, and the Sackler Center for Computational Molecular \& Materials Science, Tel Aviv University, 6997801, Tel Aviv, Israel}}

\affiliation[2]{organization={Department of Statistics and Operations Research, School of Mathematical Sciences, Tel Aviv University, 6997801, Tel Aviv, Israel}}


%%%%%%%%%%%%%%%%%%%%%%%%%%%%%%%%%%%%%%%%%%%%

\begin{abstract}
Service time fluctuations heavily affect the performance of queueing systems, causing long waiting times and backlogs. Recently, it was shown that when service times are solely determined by the server, service resetting can mitigate the deleterious effects of service time fluctuations and drastically improve queue performance \parencite{Resettingqueues}. Yet, in many queueing systems, service times have two independent sources: the intrinsic server slowdown ($S$) and the jobs' inherent size ($X$). In these, so-called $S\&X$ queues \parencite{S&Xmodel}, service resetting results in a newly drawn server slowdown while the inherent job size remains unchanged. Remarkably, resetting can be useful even then. To show this, we develop a comprehensive theory of $S\&X$ queues with service resetting. We consider cases where the total service time is either a product or a sum of the service slowdown and the jobs' inherent size. For both cases, we derive expressions for the total service time distribution and its mean under a generic service resetting policy. Two prevalent resetting policies are discussed in more detail. We first analyze the constant-rate (Poissonian) resetting policy and derive explicit conditions under which resetting reduces the mean service time and improves queue performance. Next, we consider the sharp (deterministic) resetting policy. While results hold regardless of the arrival process, we dedicate special attention to the $S\&X$-M/G/1 queue with service resetting, and obtain the distribution of the number of jobs in the system and their sojourn time. Our analysis highlights situations where service resetting can be used as an effective tool to improve the performance of $S\&X$ queueing systems. Several examples are given to illustrate our analytical results, which are corroborated using numerical simulations.

\end{abstract}

\maketitle
\section{Introduction}

Waiting in line is a frustrating yet common experience. Whether it is when standing in a long supermarket line or when waiting for a traffic light to turn green, queues appear ubiquitously in our life. Prominent examples include: call centers \parencite{Call1, Call2}, airplane boarding \parencite{Plane1, Plane3}, telecommunication and computer systems \parencite{Tele1,Mor-book}, production and manufacturing lines \parencite{Askin-Book}, enzymatic and metabolic pathways \parencite{enzymatic3,enzymatic4,enzymatic6,enzymatic8}, gene expression \parencite{Gene1,Gene2,Gene4,Gene5,Gene6}, and transport phenomena \parencite{Transport2,Transport3,Transport7,Transport9,Transport11}. 

The waiting time in a queue heavily depends on the service rate and on the arrival rate of customers or jobs. Fluctuations in service time are another important factor that affect the waiting time. While the cashier at a supermarket works at a (roughly) constant rate, other servers, e.g., computer systems \parencite{Mor-book}, and molecular machines like enzymes \parencite{Large_fluc1, Large_fluc3, Large_fluc5}, often display pronounced fluctuations in service times. These large fluctuations are a major source of backlogs and delays \parencite{Takagi, haviv2013queues,vesilo2022scaling, gupta2006fundamental}. 

To mitigate the detrimental effects large fluctuations in service times have on queue length and waiting times, various strategies have been developed. One major class is scheduling policies \parencite{harchol1999choosing, Takagi, Scheduling1}. A scheduling policy is a rule which dictates the order in which jobs in a queue are served. For example, in a single server queue, the server can reduce the overall mean waiting time by serving shorter jobs first \parencite{SJF}. In cases where service can be stopped and resumed later on, performance can be further improved by applying the shortest remaining time discipline \parencite{SRPT, SRPT-OPT}. 

Scheduling policies have proven useful in a wide variety of scenarios, yet they suffer from clear limitations. For example, the implementation of size-based scheduling policies is impossible when jobs' sizes are unknown a priori. A prominent example is given by situations where fluctuations in service times are \emph{intrinsic} to the server itself. Such scenarios are fairly common, e.g., in stochastic optimization algorithms which can take markedly different times to run on two instances of the \textit{exact same} problem \parencite{SO}. Similarly, the loading time of an internet page is stochastic, due to the varying workload on the server, although the size of the page itself does not change \parencite{ResettingInternet}. In such cases, new strategies are needed in order to tackle the problems caused by large service time fluctuations.

Resetting is the action of terminating a process and bringing it back to its initial state. Though counter-intuitive, we use resetting almost every day: when constantly refreshing a web page until it loads, or when rebooting a server to improve its performance. One can show that resetting will always prolong the completion of processes with a fixed, i.e., deterministic, completion time \parencite{branching}. However, when resetting is applied to random processes, whose completion time varies from trial to trial, the outcome can be different from the deterministic case.

The ability to reduce the mean and fluctuations of the completion time of various stochastic processes is a hallmark of resetting. The study of stochastic resetting and its applications dates back to the 1990s and the 2000s, when it was utilized in speeding up computer algorithms \parencite{Luby} and optimizing queueing systems \parencite{OR8}. Recently, emanating from the canonical diffusion problem by Evans and Majumdar \parencite{Restart1}, resetting has been shown to expedite the completion of many stochastic processes, which would otherwise take longer time to finish \parencite{Restart2, ReuveniPRL16, PalReuveniPRL17}. This happens, for example, when the completion time distribution has a decreasing failure rate. More generally, the ability of resetting to expedite stochastic processes can be traced back to the inspection paradox \parencite{IP}.

In recent years,  stochastic resetting has become a focal point of scientific interest and the subject of vigorous studies in statistical physics \parencite{Restart1,Restart2,Restart3,Restart4,Restart5,Restart6,PalJphysA,renewal-KPZ,expt-1,expt-2}, chemical and biological physics \parencite{ReuveniEnzyme1, ReuveniEnzyme2, bio-1,bio-2,bio-3, queue-input}, operations research \parencite{OR1, OR2, OR3, OR4, OR5, OR6, OR7}, reliability theory \parencite{reliability1,reliability2}, economics \parencite{income1, income2, income3} and other cross-disciplinary fields. %However, the study of stochastic resetting and its applications dates back to the '90s and early 2000s, where it was utilized in speeding up computer algorithms \parencite{Luby,algorithm}. 
For an extensive review of stochastic resetting and its applications, we refer the readers to \parencite{Review}.

In the context of queueing systems, Bonomo et al. explored the prospects of mitigating long queues and waiting times with service resetting  \parencite{Resettingqueues}. It was shown that simple service resetting mechanisms, e.g., Poissonian (resetting at a constant rate) and sharp  (resetting at constant time intervals), can reverse the deleterious effects of large fluctuations in service times. Namely, when service time fluctuations are solely \emph{intrinsic} to the server, service resetting may reduce the mean and variance of the service time and improve queue performance. Those findings imply that resetting, which does not require knowledge of job sizes, can replace size-based scheduling policies when these are not applicable.

However, like scheduling, resetting has its own limitations. A prominent one is that resetting cannot help when service times are mainly determined by factors \emph{extrinsic} to the server. A trivial example is a supermarket cashier. There, the teller serves the customers at a (roughly) constant rate, and the service duration is thus determined by the number of items collected by a customer before arriving at the counter. Evidently, resetting service and starting the checkout process anew would result in wasted time. Consequently, one should avoid resetting in cases where service times are completely determined by factors extrinsic to the server.

%This aligning with our daily experience with refreshing, namely resetting, the loading process of web pages. These resetting mechanisms are common in computing systems, for example, periodic garbage collection on servers, or annual iOS updates.

While the effect of resetting in situations where service times are solely determined either by extrinsic or by intrinsic factors is now well understood, little is known about the general case where service times have contributions from both intrinsic and extrinsic factors. In traditional queueing theory, a single random variable is assigned to represent service times, making it impossible to differentiate between extrinsic and intrinsic components. This issue was highlighted and addressed by Gardner et al. \parencite{S&Xmodel} who introduced the $S\&X$ queueing model where the service time is composed of contributions from two independent sources. 
The first source, the so-called server slowdown, is intrinsic to the server and denoted by $S$. The second source, is the jobs’ inherent size $X$ (hence the name of the model), which is extrinsic to the server.

Although the service time in the $S\&X$ model can be an arbitrary function of the intrinsic and extrinsic components, there are cases of particular interest. Consider, for example, downloading a file or loading a web page. There, the size of a single task the machine executes, namely $X$, may vary depending on the size of the web page or file. In addition, the rate at which the machine processes the task fluctuates between service attempts. This can be due to different intrinsic factors (e.g., utilization, memory usage) or because jobs are routed via different servers until they reach their final destination. Thus, the total time taken to process a task can be modeled as $S \times X$, where $S$ is the \emph{random} service time per unit of work. Note that when service resetting is applied to such queueing systems (e.g., by refreshing the browser) it will work differently on the two components of the service time. While service resetting does not change the job at hand or its size, it may change the rate at which it is served as this rate fluctuates stochastically between one service trial to another.

%%%%%%%%%%%%%%%%%%%%%%%%%%%%%%%%%%%%%%%%%%

\begin{figure}[t!]
\centering
\includegraphics[width=0.75\textwidth]{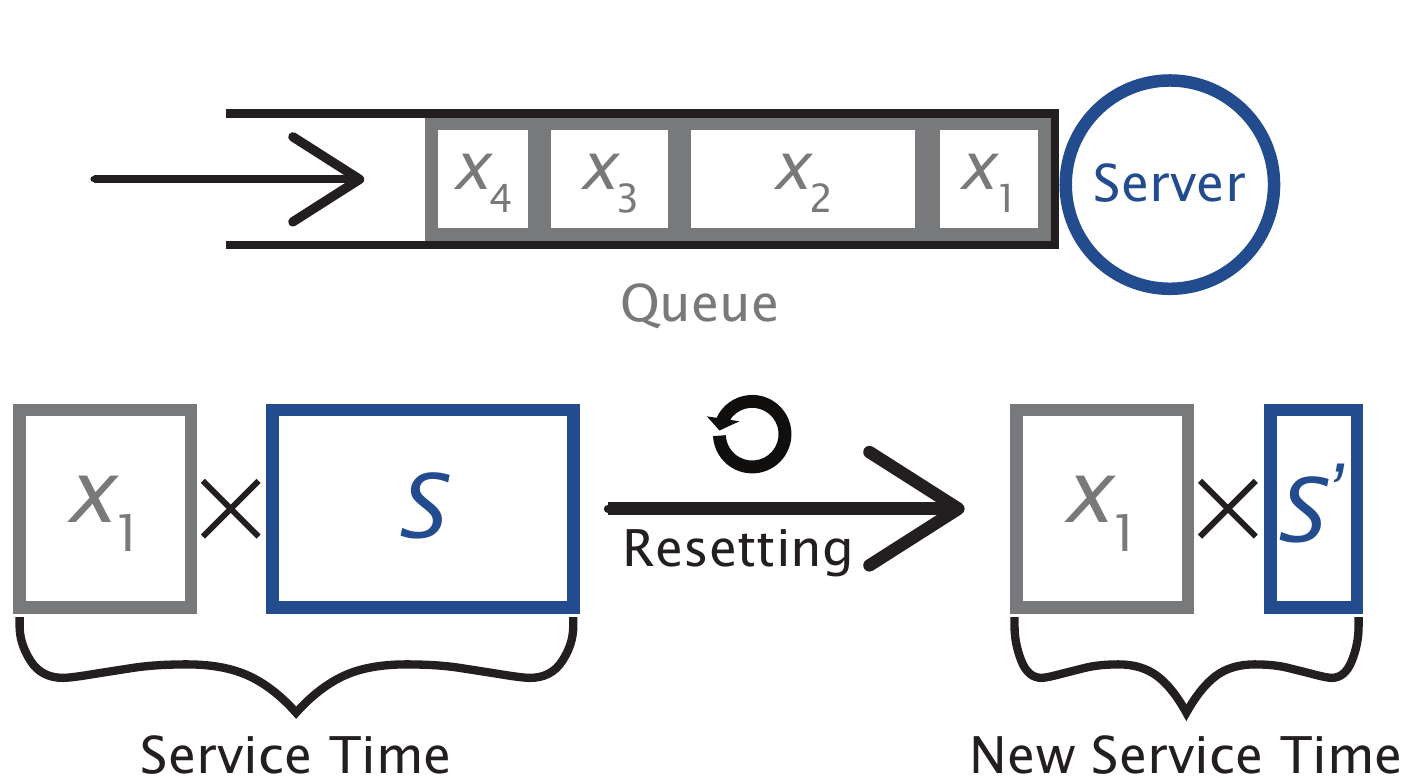}
\caption{$S\&X$ model with service resetting. Jobs of sizes $x_1, x_2, x_3, ...$ arrive at a service station. The initial required service time of the $i$-th job, $f(x_i,S)$, is a  function of the inherent size $x_i$ and the service slowdown $S$. Here, we illustrate the multiplicative case where $f(x_i,S) = x_i \times S$. A service resetting policy is applied using a timer $R$ that can be deterministic or random. If $f(x_i,S) < R$ the service is completed at the required service time. Otherwise, service is reset and a new server slowdown $S'$ is drawn. This results in a new service time $f(x_i,S')$ that is paired with a newly drawn resetting time $R'$. This process is repeated until service is completed.}
\label{Fig1a}
\end{figure}

%%%%%%%%%%%%%%%%%%%%%%%%%%%%%%%%%%%%%%%%%%

To study the effects of service resetting on $S\&X$ queues, we consider a queueing process with a generic service time that is composed of the two sources, $S$ and $X$, as discussed above. Following Gardner et al., we let $S$ and $X$ be two independent, generic random variables,  which represent the server's slowdown and the job's size, respectively. Resetting is introduced using a random time $R$, independent of both $X$ and $S$. Namely, if service is not completed by the resetting time $R$, it is halted, any previous progress is discarded, and the service is immediately restarted. Upon resetting, a new slowdown time, $S'$, is drawn, where $S$ and $S'$ are independent and identically distributed. However, the inherent size of the job in service, $X=x$, remains the same. This combined process is repeated until service is completed (Fig. 1).

In this work, we analyze the effect of service resetting on $S\&X$ queues, generalizing results obtained in \parencite{Resettingqueues} for traditional single-server queues. The paper is structured as follows. In Sec. \ref{Sec multiplicative}, we explore a multiplicative $S\&X$ model, where the overall service time is given by a product of $S$ and $X$. For this model, we first derive the service time distribution under a generally distributed resetting time. In subsection \ref{SubSec multiplicative}, we focus on the case of Poissonian resetting, i.e., resetting at a constant rate. For this case of interest, we derive an explicit condition for resetting to be advantageous, i.e., to improve queue performance by lowering the mean service time of a job. In subsection \ref{Multi-sharp}, we shift our focus to the case of sharp resetting, i.e., resetting at constant time intervals. We demonstrate the results, namely, the reduction in the mean service time, using a web page loading example, in subsection \ref{Webpage example}. In Sec. \ref{Sec additive}, we explore an additive $S\&X$ model, where the overall service time distribution is a sum of the server slowdown $S$ and the jobs’ inherent size $X$. Here too, we derive the service time distribution under a generally distributed resetting time. In subsection \ref{SubSec additive}, we concentrate on Poissonian resetting for which we obtain an explicit condition for resetting to be beneficial. In subsection \ref{Add-sharp}, we analyze the case of sharp resetting. We emphasize that the results in Secs. \ref{Sec additive} and \ref{Sec multiplicative} hold true regardless of the arrival process. In Sec. \ref{Sec distribution}, we analyze the celebrated M/G/1 queue with service resetting, and derive the full distribution of the number of jobs in the system for both the additive and multiplicative $S\&X$-M/G/1 models. We summarize in Sec. \ref{Conclusions} where we provide conclusions and outlook. The analytical results obtained in this paper are illustrated using several examples and are further corroborated using numerical simulations.

In what follows we use $f_Z(\cdot)$, $\mathbf{E}[Z]$, Var$(Z)$, $\sigma(Z)$, and $\tilde{Z}(s)~\equiv~\mathbf{E} \left[e^{-sZ}\right]$ to denote, respectively, the probability density function (PDF), expectation, variance, standard deviation, and Laplace transform of a non-negative random variable $Z$.

\section{Multiplicative $S\&X$ queues with service resetting}
\label{Sec multiplicative}
%%%%%%%%%%%%%%%%%%%%%%%%%%%%%%%%%%%%%%%%%%

Assume that the overall service time in the absence of resetting, $V$, is given by a product of the two independent sources $S$ and $X$ that were defined above. Namely, $V = S \times X$.  Let $V_R$ denote the total service time until completion under a generic service resetting policy. Given an inherent job size $x$, the conditional service time $V_R (x) \equiv\{V_{R}|X=x\}$, obeys the following renewal equation  
\begin{equation}
V_R (x)=\begin{cases}
xS & \textit{if} \quad xS<R\\
R+V_R '(x) & \textit{if} \quad  xS\geq R,
\end{cases}
\label{Renewal-Multiplication}
\end{equation}
where $R$ is a random resetting time drawn from a distribution with density $f_R(\cdot)$, and $V_{R}'(x)$ is an independent and identically distributed copy of $V_{R}(x)$. Under the multiplicative model (indicated below by the $\times$ sign), it is useful to define the following conditional random variables
\begin{align}
S_{\times}(x) &=\{S|xS<R\}, \label{XintminX} \\
R_{\times}(x) &=\{R| xS \geq R\}, \label{RminX}
\end{align}
and their corresponding Laplace transforms $\tilde{S}_{\times}(x;s)$ and $\tilde{R}_{\times}(x;s)$ 
\begin{align}
\tilde{S}_{\times}(x;s)&=\int_{0}^{\infty}e^{-st}f_{S_{\times}}(t)\,dt, \label{XintminX LT prod} \\
\tilde{R}_{\times}(x;s)&=\int_{0}^{\infty}e^{-st}f_{R_{\times}}(t)\,dt . \label{RminX LT prod}
\end{align}
In the above, we used $f_{S_{\times}}(t)$ and $f_{R_{\times}}(t)$ as the probability density functions of the conditional random variables, respectively. These are given by
\begin{align}
&f_{S_{\times}}(t)=\frac{f_{S}(t)}{\text{Pr}(xS<R)}\int_{xt}^{\infty}f_{R}(\tau)\,d\tau, \label{XintminX PDF prod} \\
&f_{R_{\times}}(t)=\frac{f_{R}(t)}{\text{Pr}(xS\geq R)}\int_{\frac{t}{x}}^{\infty}f_{S}(\tau)\,d\tau. \label{RminX PDF prod}
\end{align}
Taking expectation of both sides of Eq. (\ref{Renewal-Multiplication}) and rearranging terms, one obtains the mean conditional service time under resetting
\begin{equation}
\mathbf{E}[V_R (x)] =
x\mathbf{E}[S_{\times}(x)] +\frac{ 1 - \text{Pr}(xS< R)}{\text{Pr}(xS< R)}\mathbf{E}[R_{\times}(x)], \label{mean conditional multiplicative}
\end{equation}
where $\mathbf{E}[S_{\times}(x)]$ and $\mathbf{E}[R_{\times}(x)]$ are the means of the random variables $S_{\times}(x)$ and $R_{\times}(x)$,  respectively.

Integrating Eq. (\ref{mean conditional multiplicative}) with respect to the probability density of the jobs' inherent size, $f_{X}(x)$, one obtains the mean service time under resetting
\begin{align}
\mathbf{E}[V_R] = \int_{0}^{\infty} x  \mathbf{E}[S_{\times}(x)]~f_X(x) \, dx  +\int_{0}^{\infty} \frac{ \text{Pr}(xS\geq R)}{\text{Pr}(xS< R)} \mathbf{E}[R_{\times}(x)]~f_X(x) \, dx . \label{General Mean X}
\end{align}
The equation above for the mean service time under resetting can be explained as follows. Each service duration under resetting is composed of two contributions. The first term on the right-hand side of Eq. (\ref{General Mean X})  accounts for the mean service time of the last service trial. The second term accounts for the time wasted in failed service attempts. Namely, this is the total service time of all uncompleted service trials.

We now turn to derive the Laplace transform of $V_R$. Taking the Laplace transform of both sides of Eq.~(\ref{Renewal-Multiplication}) and rearranging terms we obtain
\begin{equation}
\tilde{V}_{R}(x;s)=\frac{\tilde{S}_{\times}(x;sx)\text{Pr}(xS<R)}{1-\tilde{R}_{\times}(x;s)\text{Pr}(xS\geq R)}. \label{LTSRa-Multiplication1}
\end{equation}

\noindent Integrating Eq. (\ref{LTSRa-Multiplication1}) with respect to the probability density of the jobs' inherent size, $f_{X}(x)$, one obtains the following Laplace transform of the service time under generic resetting in the multiplicative model
\begin{equation}
\tilde{V}_R(s)=\int_{0}^{\infty}\frac{\tilde{S}_{\times}(x;sx)\text{Pr}(xS<R)}{1-\tilde{R}_{\times}(x;s)\text{Pr}(xS\geq R)}\,f_{X}(x)\,dx . \label{Main-Multiplication}
\end{equation}

\subsection{Constant-rate Poisson resetting}
\label{SubSec multiplicative}
Consider the above-detailed service process, and further assume that inter-resetting times are taken from an exponential distribution with rate parameter $r$. Namely, the PDF of the resetting time is given by
\begin{equation}
f_R(t) = re^{-rt}, \label{Exponential PDF}
\end{equation}
for $t\geq0$. In this case, we have
\begin{align}
\text{Pr}(xS<R) = \int_{0}^{\infty} f_R(t) \text{Pr}(xS<t) \,dt = \int_{0}^{\infty} re^{-rt} \left( \int_{0}^{\frac{t}{x}} f_S(\tau) \,d\tau \right) \,dt = \tilde{S}(rx), \label{Prob X}
\end{align}
where $\tilde{S}(rx)$ is the Laplace transform of $S$ evaluated at $rx$. The last step in Eq. (\ref{Prob X}) follows by changing the order of integration.

Substituting the above result into Eqs. (\ref{XintminX PDF prod}) and (\ref{RminX PDF prod}), allows us to obtain neat expressions for the Laplace transforms of the conditional random variables $S_{\times}(x)$ and $R_{\times}(x)$. These are given by
\begin{align}
\tilde{S}_{\times}(x;s) &= \int_{0}^{\infty}e^{-st}f_{S_{\times}}(t)\,dt = \frac{\int_{0}^{\infty}e^{-st} f_{S}(t) \left(\int_{xt}^{\infty}re^{-r\tau}\,d\tau \right) \,dt}{\text{Pr}(xS<R)} \nonumber \\
&= \frac{\int_{0}^{\infty}e^{-t(s+rx)} f_{S}(t) \,dt}{\tilde{S}(rx)} = \frac{\tilde{S}(s+rx)}{\tilde{S}(rx)}, \label{LT XintminX}
\end{align}
\begin{align}
\tilde{R}_{\times}(x;s) &= \int_{0}^{\infty}e^{-st}f_{R_{\times}}(t)\,dt = \frac{\int_{0}^{\infty}e^{-st}~re^{-rt} \left(\int_{\frac{t}{x}}^{\infty}f_{S}(\tau)\,d\tau \right) \, dt}{\text{Pr}(xS\geq R)} \nonumber \\
&= \frac{r\int_{0}^{\infty}e^{-t(s+r)} \left(\int_{\frac{t}{x}}^{\infty}f_{S}(\tau)\,d\tau \right) \, dt}{1-\tilde{S}(rx)} = \frac{r\int_{0}^{\infty}e^{-t(s+r)} \left(1 - \int_{0}^{\frac{t}{x}}f_{S}(\tau)\,d\tau \right) \, dt}{1-\tilde{S}(rx)} \nonumber \\ 
&= \frac{r}{r+s}\frac{1-\tilde{S}((s+r)x)}{1-\tilde{S}(rx)}, \label{LT RminX}
\end{align}
where in the last step in Eq. (\ref{LT RminX}) we once again used the known formula for the Laplace transform of a time-domain integration: $\int_0^{\infty} dt~ \left (\int_0^{t}d\tau~g(\tau) \right ) e^{-rt} = \frac{\tilde{g}(r)}{r}$, where $\tilde{g}(r)$ denotes the Laplace transform of $g(t)$.

Substituting Eqs. (\ref{Prob X}-\ref{LT RminX}) into Eq. (\ref{Main-Multiplication}) yields
\begin{equation}
\tilde{V}_r(s)=\int_{0}^{\infty}\frac{(s+r)\tilde{S}((s+r)x)}{s+r\tilde{S}((s+r)x)}\,f_{X}(x)\,dx . \label{Main-Multiplication-Exp}
\end{equation}
where, with a slight change of notation, $\tilde{V}_r(s)$ is the Laplace transform of the service time under a constant resetting rate $r$. The moments of the service time under resetting can now be readily computed from Eq. (\ref{Main-Multiplication-Exp}), by noting that 
\begin{equation}
\mathbf{E}[V_r^n]=(-1)^n \frac{d^n}{ds^n} \tilde{V}_r(s)|_{s \to 0}. \label{Moments LT}
\end{equation}
For the first moment, we obtain
\begin{equation}
\mathbf{E}[V_r] =\int_{0}^{\infty}\frac{1-\tilde{S}(rx)}{r\tilde{S}(rx)}f_{X}(x)\,dx. \label{Mean-Multiplication-Exp}
\end{equation}

To better understand the effect of resetting on the mean service time, we consider the introduction of an infinitesimal resetting rate $\delta r$. Expanding Eq. (\ref{Mean-Multiplication-Exp}) in a Taylor series to first order in $\delta r$, we obtain (Appendix \ref{Appendix B}) 
\begin{align}
\mathbf{E} [V_{\delta r}] = \mathbf{E}[X] \mathbf{E}[S] +\delta r\mathbf{E}[X^2]\Big( \mathbf{E}[S]^{2} -\frac{1}{2}\mathbf{E}[S^2]\Big) +O(\delta r^2). \label{mean service time expansion X}
\end{align}
The first term on the right hand side of Eq. (\ref{mean service time expansion X}) is the mean of the \textit{original} service time, i.e., without resetting. The second term gives the first order correction. Rearranging the second term, we find that the introduction of a small resetting rate reduces the mean service time when the contribution of this term is negative. That is, when
\begin{equation}
CV_{S}>1, \label{CV condition X}
\end{equation}
where $CV_{S}$ is the coefficient of variation of the random variable $S$
\begin{equation}
CV_S=\sigma(S)/\mathbf{E}[S]. \label{CV}
\end{equation}

The condition in Eq. (\ref{CV condition X}) is satisfied when $\sigma(S)$, the standard deviation of the service time's intrinsic component in the absence of resetting, is larger than $\mathbf{E}[S]$, the mean of the intrinsic component. For example, if $S$ is drawn from a distribution with a decreasing failure rate, the condition in Eq. (\ref{CV condition X}) will always hold \parencite{Barlow-book}. In such cases, the introduction of a small resetting rate is guaranteed to lower the mean service time. Note that the resetting condition in Eq. (\ref{CV condition X}) depends solely on the service slowdown $S$ and is independent of the jobs' inherent size $X$. We stress that this is a property of multiplicative $S\&X$ queues, and that it does not hold in general as we show in the next section.

\begin{figure*}[t!]
\centering
\includegraphics[width=1\textwidth]{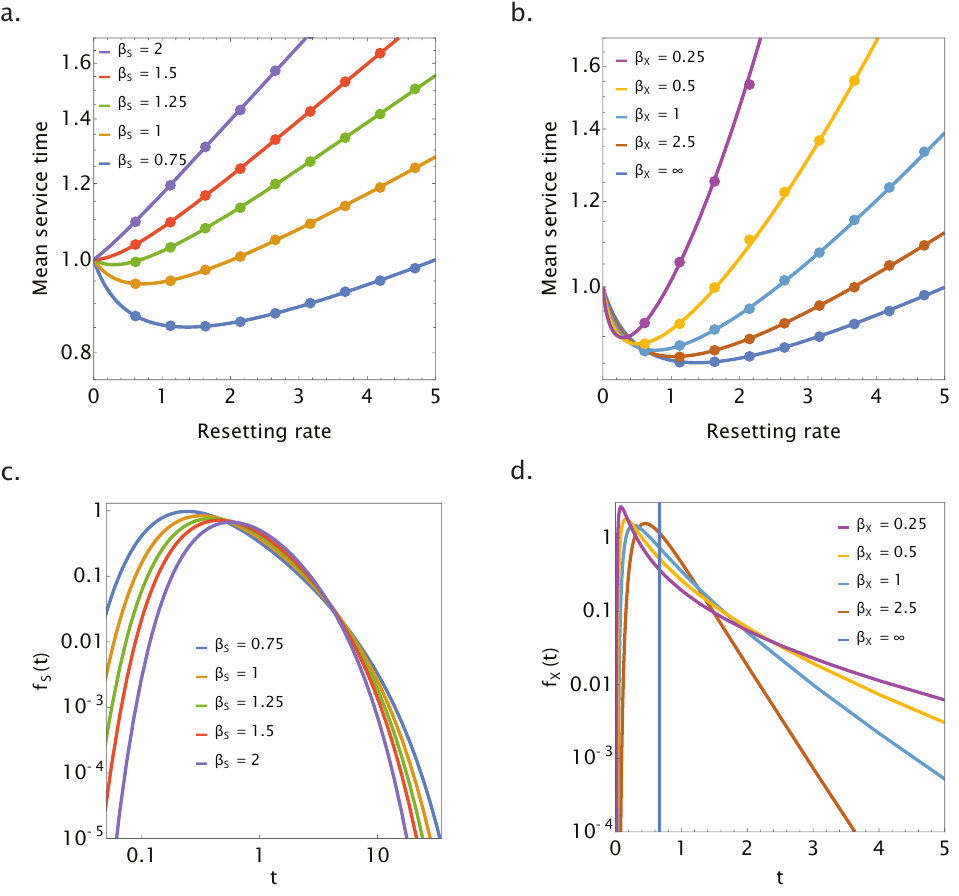}
\caption{The effect of constant-rate resetting in the multiplicative model. Panel (a): The mean service time as a function of the resetting rate. Plots are made using Eq. (\ref{Mean-Multiplication-Exp}) for a fixed inherent job size $X=x=2/3$ and a server slowdown $S$ that is taken from the Inverse Gaussian distribution ($IG$) with density $f_S(t) = \sqrt{\beta _S / 2\pi t^3} e^{-\beta _S (t-\alpha _S)^2 / 2\alpha _S ^2 t},~t>0$. Here, we set  $\mathbf{E}[S] = \alpha_S=3/2$ and vary $\beta _S$ to control the variance $\sigma^2(S)=\alpha_S ^3/\beta _S= 27\beta _S ^{-1}/8$, as illustrated in panel (c). Observe that the introduction of a small resetting rate lowers the mean service time only when the condition in Eq. (\ref{CV condition X}) is met, i.e., for $\beta_S < 3/2$. In such cases, the service time is minimized at some intermediate resetting rate. Analytical results (in solid lines) are corroborated with numerical simulations (full circles). Panel (b): The mean service time as a function of the resetting rate when both $S$ and $X$ come from the Inverse Gaussian distribution. Here, as in panel (a), we take $S\sim IG(\alpha_S = 3/2, \beta_S)$ and further set $\beta_S = 3/4$. In addition, we take $X\sim IG(\alpha_X, \beta_X)$ such that $\mathbf{E}[X]= \alpha_X = 2/3$ is fixed, and then tune $\beta_S$ to control the variance $\sigma^2(X)=8\beta_X ^{-1}/27$ as illustrated in panel (d). Observe that the introduction of  resetting lowers the mean service time regardless of the value $\beta_X$, as predicted by Eq. (\ref{CV condition X}).}\label{Fig4}
\end{figure*}

The condition in Eq. (\ref{CV condition X}) is illustrated in Fig. \ref{Fig4}a where we fix $x=2/3$ and plot the mean service time as a function of the resetting rate for intrinsic service times that are taken from the Inverse Gaussian distribution. Parameters are taken such that $\mathbf{E}[S]=3/2$ is fixed and $\sigma^2(S)$ is varied in magnitude to illustrate the transition predicted by Eq. (\ref{CV condition X}). Note that for $\sigma(S)>3/2$ the mean service time attains a minimum at an intermediate (optimal) resetting rate. Conversely, for lower values of $\sigma(S)$, the mean service time is monotonically increasing with the resetting rate.

Next, we illustrate that fluctuations in the intrinsic job size $X$ have no effect on whether resetting is beneficial or not in the multiplicative case. To this end, we plot in Fig. \ref{Fig4}b the mean service time as a function of the resetting rate for intrinsic job sizes of equal mean and different variance. Specifically, we once again take the distribution of the server slowdown as Inverse Gaussian with $\mathbf{E}[S]=3/2$ and $\sigma^2(S)=9/2$ (i.e., identical to the $\beta_S=3/4$ case from Fig. \ref{Fig4}a). Fixing this distribution for the server slowdown, we consider an intrinsic job size $X$ that comes from an Inverse Gaussian distribution with mean $\mathbf{E}[X]=2/3$ and a variance whose magnitude we gradually increase from an initial value of zero, which corresponds to the deterministic $X$ considered in Fig. \ref{Fig4}a. As predicted by Eq. (\ref{CV condition X}), since $CV_{S}^{2}=2$, the introduction of resetting reduces the mean service time regardless of the variance in the extrinsic job size.

%%%%%%%%%%%%%%%%%%%%%%%%%%%%%%%%%%%%%%%%%%

%%%%%%%%%%%%%%%%%%%%%%%%%%%%%%%%%%%%%%%%%%

\begin{figure}[t!]
\includegraphics[width=0.5\textwidth]{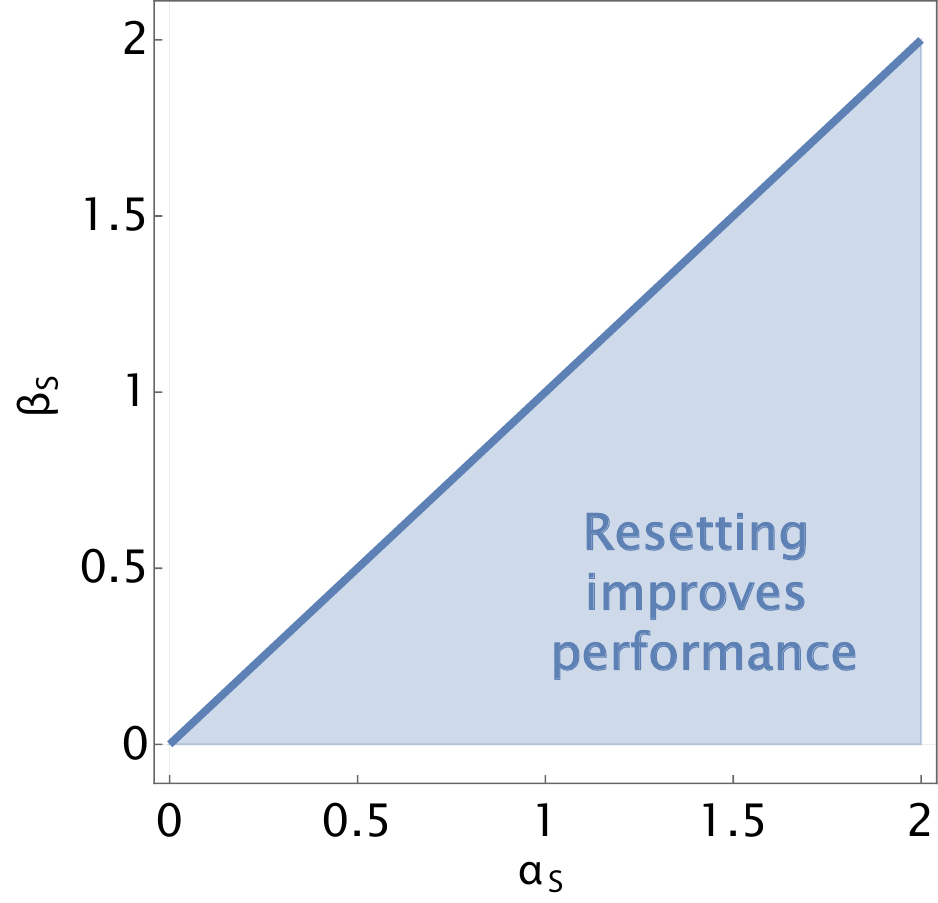}
\centering
\caption{Visualization of the condition in  Eq. (\ref{CV condition X}) for the case of Inverse Gaussian distributed $S\&X$ service times in the multiplicative model. Here, as in Fig. \ref{Fig4}, we take $S\sim IG(\alpha_S,\beta_S)$, $X\sim IG(\alpha_X,\beta_X)$. Tuning $\alpha_S$ and $\beta_S$ we control the mean $\mathbf{E}[S] = \alpha_S$, and the variance $\sigma^2(S)=\alpha_S ^3/\beta _S$, of the intrinsic component of the service time. The condition in Eq. (\ref{CV condition X}) is satisfied when $CV_{S}^{2} = \alpha_S / \beta_S$ is larger than 1. This region of the phase space, where the introduction of service resetting reduces the mean service time, is shaded in blue.}
\label{Fig5}
\end{figure}

Summarizing, we see that service resetting is beneficial when $\sigma(S)$, the standard deviation in the intrinsic component of the service time, is larger than $\mathbf{E}[S]$, the mean of the intrinsic component of the service time. The condition for resetting being beneficial is agnostic to the extrinsic component of the service time distribution. A phase space visualization of the condition in  Eq. (\ref{CV condition X}), for the case where both $S$ and $X$ are drawn from an Inverse Gaussian distribution, is provided in Fig. {\ref{Fig5}}.

%%%%%%%%%%%%%%%%%%%%%%%%%%%%%%%%%%%%%%%%%

%%%%%%%%%%%%%%%%%%%%%%%%%%%%%%%%%%%%%%%%%%

\subsection{Sharp resetting}
\label{Multi-sharp}
Consider now the case of sharp, i.e. deterministic, resetting. When the resetting time $R$ is deterministic, its probability density function is given by 
\begin{align}
    f_R(t)=\delta(t-\tau)~, \label{sharp density}
\end{align}
where $\tau$ is the resetting time and $\delta(\cdot)$ is the delta function.
In the case of sharp resetting, the PDFs of the conditional random variables $S_{\times}(x)$ and $R_{\times}(x)$ can be written as follows
\begin{align}
&f_{S_{\times}}(t)=\frac{f_{S}(t)}{\int_{0}^{\frac{\tau}{x}}f_{S}(t')\,dt'}\theta(\tau-tx), \label{SX sharp prod} \\
&f_{R_{\times}}(t)=\frac{\delta(t-\tau)}{1 - \int_{0}^{\frac{\tau}{x}}f_{S}(t')\,dt'}\int_{\frac{t}{x}}^{\infty}f_{S}(t')\,dt'. \label{RX sharp prod}
\end{align}
where we have used $\text{Pr}(R > t)=\int_t^{\infty} dt'~\delta(t'-\tau) = \theta(\tau-t)$,  where $\theta$ is the Heaviside step function. The expectations of the above random variables are then given by
\begin{align}
&\mathbf{E}[S_{\times}(x)] = \int_0^{\infty} t f_{S_{\times}}(t) dt = \frac{\int_{0}^{\frac{\tau}{x}} t f_{S}(t) dt}{\int_{0}^{\frac{\tau}{x}}f_{S}(t')\,dt'}, \label{SX sharp prod mean} \\
&\mathbf{E}[R_{\times}(x)]=\int_0^{\infty} t f_{R_{\times}}(t) dt = \frac{\tau \left( 1 - \int_{0}^{\frac{\tau}{x}}f_{S}(t')\,dt' \right)}{1 - \int_{0}^{\frac{\tau}{x}}f_{S}(t')\,dt'} = \tau. \label{RX sharp prod mean}
\end{align}
Substituting the expectation for $\mathbf{E}[S_{\times}(x)]$ and $\mathbf{E}[R_{\times}(x)]$ given above into Eq. (\ref{General Mean X}) yields the following formula for $\mathbf{E}[V_{\tau}]$, the mean service time under sharp resetting
\begin{align}
\mathbf{E}[V_{\tau}] & = \int_{0}^{\infty} \frac{\int_{0}^{\frac{\tau}{x}} t f_{S}(t) dt}{\int_{0}^{\frac{\tau}{x}}f_{S}(t')\,dt'}~x f_X(x) \, dx  +\int_{0}^{\infty} \frac{ 1- \int_{0}^{\frac{\tau}{x}}f_{S}(t)\,dt}{\int_{0}^{\frac{\tau}{x}}f_{S}(t')\,dt'} \tau~f_X(x) \, dx \nonumber \\ 
& = \int_{0}^{\infty} \frac{\int_{0}^{\frac{\tau}{x}} t f_{S}(t) dt}{\int_{0}^{\frac{\tau}{x}}f_{S}(t')\,dt'}~x f_X(x) \, dx  +\tau \left( \int_{0}^{\infty} \frac{ f_X(x)}{\int_{0}^{\frac{\tau}{x}}f_{S}(t')\,dt'} \, dx - 1 \right ). \label{Sharp Mean X}
\end{align}

\subsection{Web page loading time example}
\label{Webpage example}
As previously discussed, the multiplicative model applies to queueing systems where each job requires $X$ units of work taken from a distribution $f_X(x)$. The inherent size of the job in service, $X = x$, is fixed throughout the sojourn of the job in the queue. Simultaneously, each unit of work requires $S$ units of time for completion. This time is intrinsic to the server and may vary between service attempts, or alternatively, when a job is directed to a statistically identical server. Hence, the time needed to serve a generic job, denoted as $V$, can be calculated as $S \times X$.

These scenarios are common in computing systems, such as when downloading a file or loading a web page. There, the size of a single task executed by the machine is fixed. However, the rate at which the machine runs the task fluctuates due to different intrinsic factors (utilization, memory usage, etc.). Thus, the time taken to process the same task varies between different service cycles. Next, we use the model of multiplicative $S\&X$ queues to demonstrate how resetting can expedite the execution of tasks in these scenarios.

\begin{figure}[t!]
\includegraphics[width=1\textwidth]{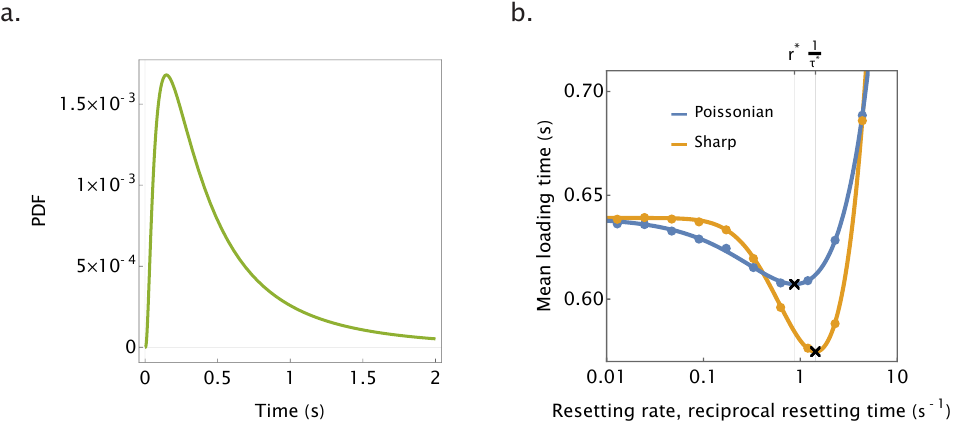}
\centering
\caption{Loading time in seconds of a web page under Poissonian and sharp 'refreshing'. Panel (a): The loading time distribution of a single, randomly chosen web page, is given by a log-normal distribution with parameters $\mu = 5.97$ and $\sigma = 0.99$ \parencite{ResettingInternet}. Panel (b): The mean loading time with Poissonian (blue) and sharp (orange) resetting (browser refresh), as a function of the resetting rate, or reciprocal of the sharp resetting time. Plots are made using Eq. (\ref{Mean-Multiplication-Exp}) and Eq. (\ref{Sharp Mean X}), for a fixed inherent job size $x= 0.9\cdot \mathbf{E}[T] \simeq 575.184$ and a server slowdown $S$ that is taken from the log-normal distribution, with parameters $\mu = 5.97 - \ln(575.184) \simeq -0.385$, and $\sigma = 0.99$. The mean loading times at the optimal resetting rate, $r^*$, and the reciprocal of the optimal resetting time, $\tau^*$, are indicated by $\mathbf{x}$ symbols. Analytical results in solid lines are corroborated with numerical simulations (full circles).}
\label{WebPageEx}
\end{figure}

We consider the scenario of a computer system processing randomly chosen web pages in succession. In \parencite{ResettingInternet}, it is shown that the download time of the index.html file on the main page of different websites, which corresponds to their loading time, follows a log-normal distribution with parameters $\mu = 5.97$ and $\sigma = 0.99$ (Fig. \ref{WebPageEx}a). In general, the PDF of a log-normally distributed random variable is given by
\begin{equation}
f_S(t)=\frac{1}{\sqrt{2\pi} \sigma t}~e^{-\frac{(\ln t-\mu)^2}{2\sigma^2}}~, \label{log normal dist}
\end{equation}
for $t>0$, where $\mu \in (-\infty, \infty)$ and $\sigma>0$.

Given the above result, we set the loading time $T$, of a single, randomly chosen web page, to be taken from a log-normal distribution with parameters $\mu = 5.97$ and $\sigma = 0.99$. For illustration purposes, we take a fixed inherent job size $X=x$. Assuming a multiplicative model scenario, the service slowdown, $S$, is given by $S=T/x$. It follows that $S$ is log-normally distributed, with parameters $\mu = 5.97 - \ln(x)$ and $\sigma = 0.99$.

%%%%%%%%%%%%%%%%%%%%%%%%%%%%%%%%%%%%%%%%%%

%%%%%%%%%%%%%%%%%%%%%%%%%%%%%%%%%%%%%%%%%%%%%%%%%

We now introduce stochastic resetting into the service, namely, `refreshing' the browser after some random or deterministic time. For constant-rate Poissonian resetting, the mean service time is given by Eq. (\ref{Mean-Multiplication-Exp}), which requires the Laplace transform of the log-normal distribution. The latter, does not have an analytical closed-form, but it can be evaluated numerically for any choice of parameters.

In the case of sharp resetting, the mean service time is given by  Eq. (\ref{Sharp Mean X}), which requires the PDFs of both the extrinsic and intrinsic service times. Here, the former is given by $\delta(t-x)$, and the latter is given by Eq. (\ref{log normal dist}). With these PDFs at hand, the mean service time under resetting can be computed by numerical evaluation of the required integrals, similarly to the constant-rate resetting case.

For example, in Fig. \ref{WebPageEx}b, we set $x = 0.9\cdot \mathbf{E}[T] \simeq 575.184$ and plot the mean service time under Poissonian and sharp resetting. In both cases, a minimum is attained at an optimal resetting rate or time, depending on the resetting scheme. Observe that the optimal mean service time under sharp resetting is lower than that obtained for Poissonian resetting. Interestingly, in both cases we find an optimal refreshing rate of about $1 s$, which aligns with common everyday experience of surfing the web.

\section{Additive $S\&X$ queues with service resetting}
\label{Sec additive}

%%%%%%%%%%%%%%%%%%%%%%%%%%%%%%%%%%%%%%%%%%

\begin{figure}[t!]
\centering
\includegraphics[width=0.75\textwidth]{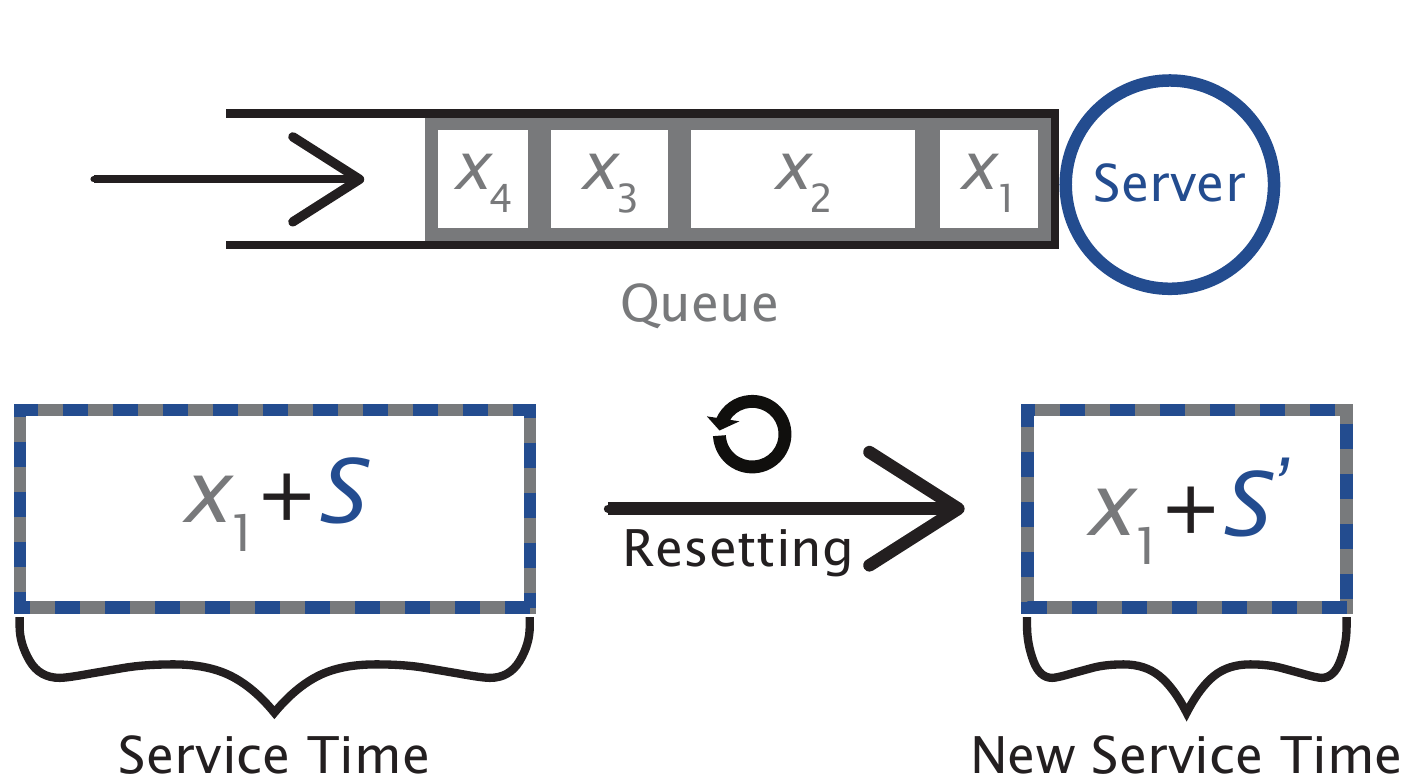}
\caption{The additive $S\&X$ model with service resetting. Jobs of sizes $x_1, x_2, x_3, ...$ arrive at a service station. The initial required service time of the $i$-th job is given by  $f(x_i,S) = x_i + S$. A service resetting policy is applied using a timer $R$ that can be deterministic or random. If $f(x_i,S) < R$ the service is completed at the required service time. Otherwise, service is reset and a new server slowdown $S'$ is drawn. This results in a new service time $x_i + S'$ that is paired with a newly drawn resetting time $R'$. This process is repeated until service is completed.}
\label{Fig1b}
\end{figure}

%%%%%%%%%%%%%%%%%%%%%%%%%%%%%%%%%%%%%%%%%%

In this section, we introduce and analyze the additive $S\&X$ model with service resetting. In this model, the overall service time is a sum of the server slowdown $S$ and the jobs’ inherent size $X$. If service is not completed before the resetting time $R$, it is halted and immediately restarted, discarding any previous progress. The new overall service time is then given by $S'+X$, where $S'$ is an independent and identically distributed copy of $S$, while the inherent size of the job in service, $X=x$, remains unchanged (see Fig. \ref{Fig1b}).

The description above implicitly assumes that the service of a job is carried out in a single step of $S+X$ units of time. Namely, we assume that service cannot be broken into steps where the $X$ component is served first and the $S$ component is served second (or vice versa). We note that the latter scenario can be mapped onto the model introduced and analyzed in Bonomo et al. \parencite{Resettingqueues}. In this variant, there is no point resetting in the $X$ phase, and the work associated with the inherent job size is carried out uninterrupted. Thus, only the step associated with the service slowdown is restarted, which gives a total service time of $X+S_R$, where $S_R$ is the service slowdown time under resetting. We thus consider this case solved and focus on the case where the total service time consists of a single phase that cannot be broken down into its individual components.

Scenarios where the intrinsic and extrinsic service components are intertwined and cannot be disentangled are commonly encountered. Consider, for example, a case where on top of the intrinsic service time of a job one also needs to wait a random time $S$ during which service is halted. For example, this could be because the server requires rest periods (breaks), is called to attend more urgent tasks, or is simply malfunctioning and requires repair at different time points along the service process. The customer may then choose to stop service and start it over with a different server that is available. Restarting could also make sense when the time $S$ is attributed to some portion of the job that is delegated by the server to a third-party server. As some third-party servers are faster than others, it may be worthwhile to switch between them from time to time --- even if this means that all work done so far must be discarded.

Motivated by the above, we denote the overall service time in the absence of resetting as, $B = S + X$. Similar to the treatment of the multiplicative model, we let $B_R$ denote the total service time under resetting. Given an inherent job size of $x$, the conditional service time $B_{R}(x)\equiv\{B_{R}|X=x\}$ obeys the following renewal equation
\begin{equation} 
B_{R}(x)=\begin{cases}
x+S & \textit{if} \quad x+S<R\\
R+B_{R}'(x) & \textit{if} \quad x+S\geq R,
\end{cases}
\label{Renewal-Addition}
\end{equation}
\noindent
where, once again, $R$ is a random resetting time drawn from a distribution with density $f_R(\cdot)$, and $B_R '(x)$ is an independent and identically distributed copy of $B_R (x)$. This combined process, which is illustrated in Fig. \ref{Fig1b}, is repeated until service is completed. Under the additive model (indicated by the $+$ sign), it is useful to define the following conditional random variables
\begin{align}
S_{+}(x)&=\{S|x+S<R\}, \label{Xintmin+} \\
R_{+}(x)&=\{R|x+S\geq R\}, \label{Rmin+}
\end{align}
and their corresponding Laplace transforms $\tilde{S}_{+}(x;s)$ and $\tilde{R}_{+}(x;s)$
\begin{align}
\tilde{S}_{+}(x;s)&=\int_{0}^{\infty}e^{-st}f_{S_{+}}(t)\,dt , \label{XintminX LT} \\
\tilde{R}_{+}(x;s)&=\int_{0}^{\infty}e^{-st}f_{R_{+}}(t)\,dt . \label{RminX LT}
\end{align}
The PDFs of these conditional random variables are given by
\begin{align}
&f_{S_{+}}(t)=\frac{f_{S}(t)}{\text{Pr}(x+S<R)}\int_{t+x}^{\infty}f_{R}(\tau)\,d\tau, \label{XintminX PDF sum} \\
&f_{R_{+}}(t)=\frac{f_{R}(t)}{\text{Pr}(x+S\geq R)}\int_{t-x}^{\infty}f_{S}(\tau)\,d\tau. \label{RminX PDF sum}
\end{align}
Note that in Eq. (\ref{RminX PDF sum}), for $t<x$ the lower integration limit is set to $0$.

Taking the expectation of both sides of Eq. (\ref{Renewal-Addition}) we get 
\begin{align}
\mathbf{E}[B_{R}(x)] = ~&~\mathbf{E}[(x+S)|x+S< R]\text{Pr}(x+S< R) \nonumber \\ &+ \mathbf{E}[(R+B_{R}'(x))|x+S\geq R]\text{Pr}(x+S\geq R).
\label{mean conditional additive1}
\end{align}
Using Eqs. (\ref{Xintmin+}) and (\ref{Rmin+}) and rearranging terms in Eq. (\ref{mean conditional additive1}), one obtains the mean conditional service time under resetting
\begin{equation}
\mathbf{E}[B_{R}(x)] =
x + \mathbf{E}[S_{+}(x)] +
\frac{ 1 - \text{Pr}(x+S< R)}{\text{Pr}(x+S< R)}\mathbf{E}[R_{+}(x)], \label{mean conditional additive}
\end{equation}
where $\mathbf{E}[S_{+}(x)]$ and $\mathbf{E}[R_{+}(x)]$ are the means of the random variables $S_{+}(x)$ and $R_{+}(x)$ respectively.

Integrating Eq. (\ref{mean conditional additive}) with respect to the probability density of the jobs' inherent size, $f_{X}(x)$, one obtains the mean service time under resetting
\begin{align}
\mathbf{E}[B_R] = \, \mathbf{E}[X] + \int_{0}^{\infty} \mathbf{E}[S_{+}(x)] f_X(x) \, dx  +\int_{0}^{\infty} \frac{ \text{Pr}(x+S \geq R)}{1 - \text{Pr}(x+S \geq R)} \mathbf{E}[R_{+}(x)] f_X(x) \, dx. \label{General Mean +}
\end{align}
Note the similarities between Eq. (\ref{General Mean X}) for the multiplicative case and Eq. (\ref{General Mean +}) for the additive case. The first two terms on the right-hand side of Eq. (\ref{General Mean +}) account for the mean service time of the last service trial. The last term accounts for the time wasted in the geometric sum of  failed service attempts.

Finally, we turn to derive the Laplace transform of $B_R$. Taking the Laplace transform of both sides of Eq. (\ref{Renewal-Addition}) and rearranging terms we obtain
\begin{equation}
\tilde{B}_{R}(x;s)=\frac{e^{-sx}\tilde{S}_{+}(x;s)\text{Pr}(x+S<R)}{1-\tilde{R}_{+}(x;s)\text{Pr}(x+S\geq R)}.
\label{LTSRa-Addition1}
\end{equation}
Integrating Eq. (\ref{LTSRa-Addition1}) with respect to the probability density of the jobs' inherent size, $f_{X}(x)$, one obtains the Laplace transform of the service time under a generic resetting timer
\begin{equation}
\tilde{B}_R(s)=\int_{0}^{\infty}\frac{e^{-sx}\tilde{S}_{+}(x;s)\text{Pr}(x+S<R)}{1-\tilde{R}_{+}(x;s)\text{Pr}(x+S\geq R)}\,f_{X}(x)\,dx. \label{Main-Additive}
\end{equation}

\subsection{Constant-rate Poisson resetting}
\label{SubSec additive}
Consider the above-detailed service process and further assume that resetting times are taken from an exponential distribution with rate parameter $r$. Namely, for $t>0$, the probability density of the resetting time distribution is given by Eq. (\ref{Exponential PDF}). In this case, we have 
\begin{align}
\text{Pr}(x+S<R) &= \int_{0}^{\infty} f_R(t) \text{Pr}(x+S<t) \,dt = \int_{x}^{\infty} re^{-rt} \left( \int_{0}^{t-x} f_S(\tau) \,d\tau \right) \,dt = e^{-rx}\tilde{S}(r), \label{Prob +}
\end{align}
where $\tilde{S}(r)$ is the Laplace transform of $S$, evaluated at the resetting rate $r$. Note that in moving from the first to the second equality in Eq. (\ref{Prob +}), we used  $\text{Pr}(x+S<t) = 0$ for $t<x$. The last step then follows by a change of  variable and order of integration.

Substituting the above result into Eqs. (\ref{XintminX PDF sum}) and (\ref{RminX PDF sum}), allows us to obtain neat expressions for the Laplace transforms of the conditional random variables $S_{+}(x)$ and $R_{+}(x)$. These are given by
\begin{align}
\tilde{S}_{+}(x;s) &= \int_{0}^{\infty}e^{-st}f_{S_{+}}(t)\,dt = \frac{\int_{0}^{\infty}e^{-st}f_{S}(t)\left(\int_{t+x}^{\infty} re^{-r\tau} \,d\tau \right) \,dt}{\text{Pr}(x+S<R)} \nonumber \\
&= \frac{\int_{0}^{\infty}e^{-t(r+s)-rx}f_{S}(t) \,dt}{e^{-rx}\tilde{S}(r)} = \frac{\int_{0}^{\infty}e^{-t(r+s)}f_{S}(t) \,dt}{\tilde{S}(r)} = \frac{\tilde{S}(s+r)}{\tilde{S}(r)}, \label{LT Xintmin+}
\end{align}
\begin{align}
\tilde{R}_{+}(x;s) &= \int_{0}^{\infty}e^{-st}f_{R_{+}}(t)\,dt = \frac{\int_{0}^{\infty}e^{-st}~re^{-rt} \left(\int_{t-x}^{\infty}f_{S}(\tau)\,d\tau \right) \, dt}{\text{Pr}(x+S\geq R)} \nonumber \\
& = \frac{r\int_{0}^{\infty}e^{-(s+r)t} \left(1-\int_{0}^{t-x}f_{S}(\tau)\,d\tau \right) \, dt}{1-e^{-rx}\tilde{S}(r)} = \frac{r}{s+r}\frac{1-e^{-(s+r)x}\tilde{S}(s+r)}{1-e^{-rx}\tilde{S}(r)}, \label{LT Rmin+}
\end{align}
where in the last step in Eq. (\ref{LT Rmin+}) we used a known formula for the Laplace transform of a time-domain integration: $\int_0^{\infty} dt~ \left (\int_0^{t}d\tau~g(\tau) \right ) e^{-rt} = \frac{\tilde{g}(r)}{r}$, where $\tilde{g}(r)$ denotes the Laplace transform of $g(t)$.

Substituting Eqs. (\ref{Prob +}-\ref{LT Rmin+}) into Eq. (\ref{Main-Additive}) yields
\begin{equation}
\tilde{B}_r(s)=\int_{0}^{\infty}\frac{(s+r)e^{-(s+r)x}\tilde{S}(s+r)}{s+re^{-(s+r)x}\tilde{S}(s+r)}\,f_{X}(x)\,dx . \label{Main-Additive-Exp}
\end{equation}
The moments of the service time under resetting
can be computed from Eq. (\ref{Main-Additive-Exp}), by utilizing Eq. (\ref{Moments LT}). For the first moment, we obtain
\begin{equation}
\mathbf{E}[B_r]=\frac{\tilde{X}(-r)-\tilde{S}(r)}{r\tilde{S}(r)}. \label{Mean-Additive-Exp}
\end{equation}
Note that Eq. (\ref{Mean-Additive-Exp}) is well defined only when
\begin{equation}
\tilde{X}(-r) = \int_{0}^{\infty}e^{rx}f_X (x)\,dx  \label{regularity condition}
\end{equation}
is finite. This condition is automatically violated for jobs' inherent size distributions with tails that are not exponentially bounded, i.e. heavy-tailed distributions   \parencite{Heavy-tail}. In what follows we consider jobs' inherent size distributions, and resetting rates, for which $\tilde{X}(-r)$ is finite.

To better understand the effect of resetting on the mean service time, consider the introduction of an infinitesimal resetting rate $\delta r$. Expanding Eq. (\ref{Mean-Additive-Exp}) in a Taylor series to first order in $\delta r$, we obtain (Appendix \ref{Appendix A})  
\begin{align}
\mathbf{E}[B_{\delta r}] = & ~ \left (\mathbf{E}[X] +\mathbf{E}[S] \right ) +\delta r \Big( \mathbf{E}[X] \mathbf{E}[S] + \mathbf{E}[S]^{2} + \frac{1}{2}\left(\mathbf{E}[X^2] -\mathbf{E}[S^2] \right)\Big) +O(\delta r^2). \label{mean service time expansion +}
\end{align}
The first term on the right hand side of Eq. (\ref{mean service time expansion +}) is the mean of the \textit{original} service time, i.e., without resetting; and the second term gives the first order correction to this result. In particular, when this term is negative the introduction of a small resetting rate will reduce the mean service time. Rearranging the parenthesis in the second term, we find that this occurs  when 
\begin{equation}
\Big( \mathbf{E}[X] +\mathbf{E}[S] \Big)^{2}<\sigma^2(S) - \sigma^2(X). \label{CV condition +}
\end{equation}
Note that the resetting condition in the additive case depends not only on the service slowdown $S$ but also on the jobs' inherent size $X$. This stands in contrast to the condition in the additive case which depends solely on the service slowdown S and is independent of the jobs’ inherent size X, as can be seen in Eq. (\ref{CV condition X}).

To better understand the above condition consider a scenario where the slowdown time $S$ is taken from a generic distribution, and the jobs’ inherent size $X$ is fixed, i.e., $X=x$. In this scenario, $\sigma^2(X)=0$ and  the above condition reduces to
\begin{align}
CV_{S} > 1 + \frac{x}{\mathbf{E}[S]}.
\label{CV condition + sim}
\end{align}
The condition given in Eq. (\ref{CV condition + sim}) generalizes the $CV_{S} > 1$ condition that was previously derived in \parencite{Resettingqueues} for the case $x=0$.

%%%%%%%%%%%%%%%%%%%%%%%%%%%%%%%%%%%%%%%%%%

\begin{figure*}[t!]
\centering
\includegraphics[width=1\textwidth]{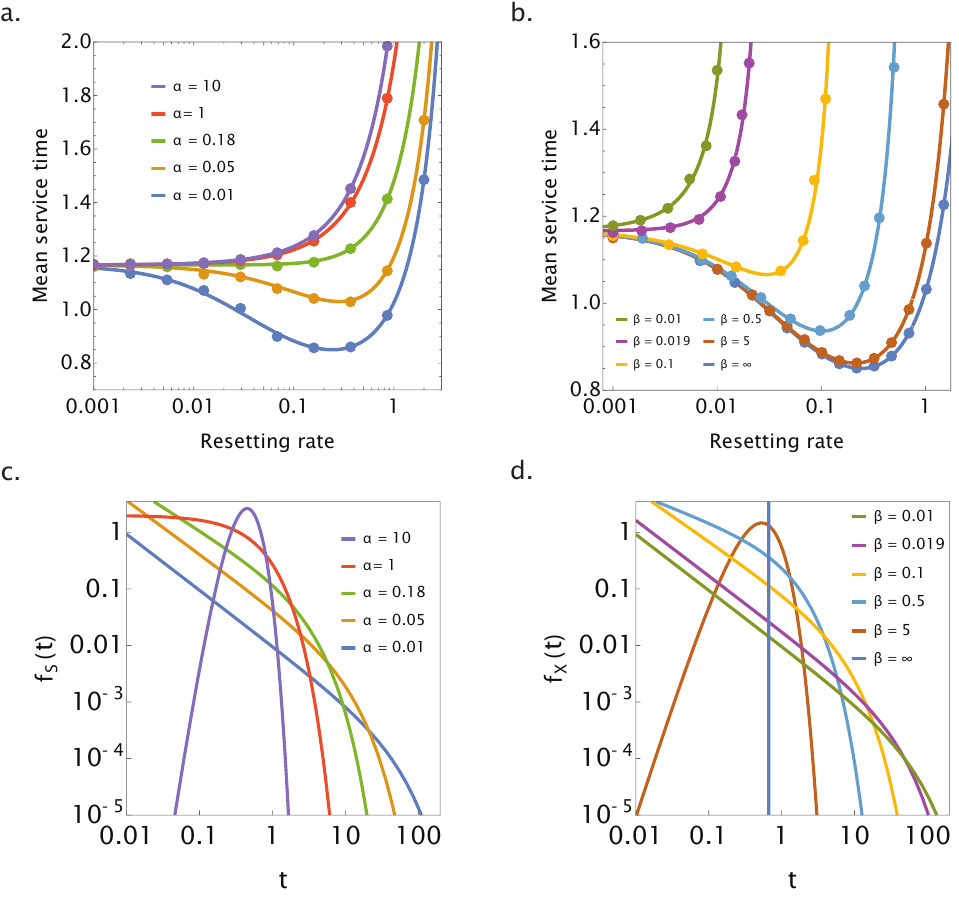}
\caption{The effect of constant-rate resetting in the additive model. Panel (a): The mean service time as a function of the resetting rate. Plots are made using Eq. (\ref{Mean-Additive-Exp}) for a fixed inherent job size $X=x=2/3$ and a server slowdown $S$ that is taken from the Gamma distribution with density $f_S(t) = \frac{1}{\Gamma(\alpha)\theta^\alpha}t^{\alpha - 1} e^{-t/\theta},~t>0$. Here, we set  $\mathbf{E}[S] = \theta\alpha=1/2$ and vary $\alpha$ to control the variance $\sigma^2(S)=\theta^2\alpha=\alpha^{-1}/4$, as illustrated in panel (c). Observe that the introduction of a small resetting rate lowers the mean service time only when the condition in Eq. (\ref{CV condition + sim}) is met, i.e., for $\alpha < 9/49 \simeq 0.184$. In such cases, the service time is minimized at some intermediate resetting rate. Analytical results (in solid lines) are corroborated with numerical simulations (full circles). Panel (b): The mean service time as a function of the resetting rate when both $S$ and $X$ come from the Gamma distribution. Here, as in panel (a), we take $S\sim\Gamma(\alpha,\theta=\frac{1}{2\alpha})$ and further set $\alpha = 0.01$. In addition, we take $X\sim\Gamma(\beta, \frac{2}{3 \beta })$ such that $\mathbf{E}[X]=2/3$ is fixed, and then tune $\beta$ to control the variance $\sigma^2(X)=4\beta^{-1}/9$ as illustrated in panel (d). Observe that the introduction of a small resetting rate lowers the mean service time only when the condition in Eq. (\ref{CV condition +}) is met, i.e., for $\beta > 16/851 \simeq 0.019$.}
\label{Fig2}
\end{figure*}

%%%%%%%%%%%%%%%%%%%%%%%%%%%%%%%%%%%%%%%%%%

Note that the condition in Eq. (\ref{CV condition + sim}) is satisfied when the standard deviation in the intrinsic component of the service time in the absence of resetting, $\sigma(S)$, is larger than the mean service time without resetting, $\mathbf{E}[B]=x +\mathbf{E}[S]$. In such cases, the introduction of a small resetting rate is guaranteed to lower the mean service time. This is illustrated in Fig. \ref{Fig2}a where we fix $x=2/3$ and plot the mean service time as a function of the resetting rate for intrinsic service times that are taken from the Gamma distribution. Parameters are taken such that $\mathbf{E}[S]=1/2$ is fixed and $\sigma^2(S)$ is varied in magnitude to illustrate the transition predicted by Eq. (\ref{CV condition + sim}). Note that for $\sigma(S)>7/6$ the mean service time attains a minimum at an intermediate (optimal) resetting rate. Conversely, for lower variances, the mean service time is monotonically increasing with the resetting rate.

Next, we illustrate the effect of fluctuations in the intrinsic job size $X$. To this end, we plot in Fig. \ref{Fig2}b the mean service time as a function of the resetting rate for intrinsic job sizes of equal mean and different variance. Specifically, we once again take the distribution of the server slowdown as Gamma with $\mathbf{E}[S]=1/2$ and $\sigma^2(S)=25$ (i.e., identical to the $\alpha=0.01$ case from Fig. \ref{Fig2}a). Fixing this distribution for the server slowdown, we consider an inherent job size $X$ that comes from a Gamma distribution with mean $\mathbf{E}[X]=2/3$ and a variance whose magnitude we gradually increase from an initial value of zero (the deterministic case considered in Fig. \ref{Fig2}a). As predicted by Eq. (\ref{CV condition +}), the introduction of resetting reduces the mean service time as long as the variance in the inherent job size is smaller that a certain threshold, which is here given by $\sigma^2(X)=23\frac{23}{36}$. In such cases, the service time is minimized at some intermediate (optimal) resetting rate. 

\begin{figure}[t!]
\includegraphics[width=0.5\textwidth]{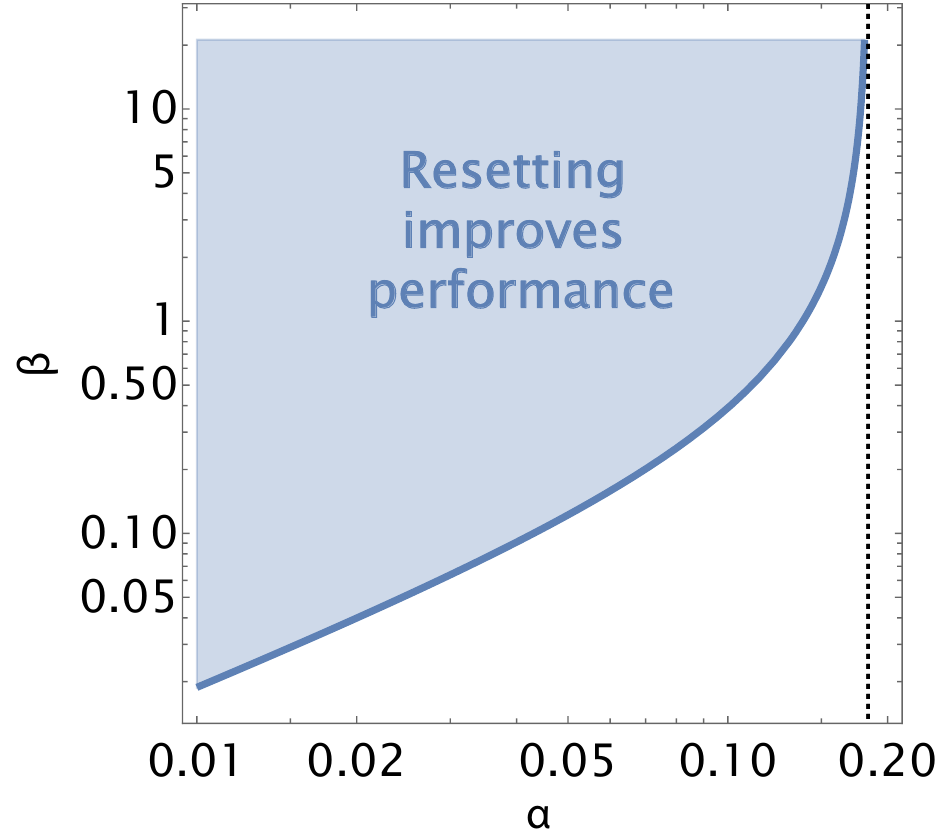}
\centering
\caption{Visualization of the condition in Eq. (\ref{CV condition +}) for the case of Gamma distributed $S\&X$ service times in the additive model. Here, as in Fig. \ref{Fig2}, we take $S\sim\Gamma(\alpha,\frac{1}{2\alpha})$ and  $X\sim\Gamma(\beta, \frac{2}{3 \beta })$ such that $\mathbf{E}[S] = 1/2$ and $\mathbf{E}[X]=2/3$ are fixed. Tuning $\alpha$ and $\beta$ we control the variance of these service time components as $\sigma^2(S)=\alpha^{-1}/4$ and $\sigma^2(X)=4\beta^{-1}/9$. The condition in Eq. (\ref{CV condition +}) is satisfied when $\sigma^2(S)$ is large enough and when $\sigma^2(X)$ is small enough. This part of the phase space, where the introduction of service resetting reduces the mean service time, is shaded in blue.}
\label{Fig3}
\end{figure}

Summarizing, we see that service resetting is beneficial when the variance in the intrinsic component of the service time is large compared to the mean-squared service time. In addition, the variance in the inherent job size must be small by the same measure. As a result, the parameter phase space is divided into two regions: one where the introduction of resetting reduces the mean service time and one where the converse happens. This phase space is illustrated in Fig. \ref{Fig3} for the case of Gamma distributed $S\&X$ service times.

%%%%%%%%%%%%%%%%%%%%%%%%%%%%%%%%%%%%%%%%%%

%%%%%%%%%%%%%%%%%%%%%%%%%%%%%%%%%%%%%%%%%%%%%%%%%

\subsection{Sharp resetting}
\label{Add-sharp}
Consider the above-detailed service process and further assume that resetting times are taken from a sharp distribution with parameter $\tau$. Namely, for $t > 0$, the probability density of the resetting time distribution is given by Eq. (\ref{sharp density}). In this case, the PDFs of the conditional random variables $S_{+}(x)$ and $R_{+}(x)$ can be written as followed
\begin{align}
&f_{S_{+}}(t)=\frac{f_{S}(t)}{\int_{0}^{\tau -x}f_{S}(t')\,dt'}\theta(\tau-(t + x)), \label{RX sharp sum} \\
&f_{R_{+}}(t)=\frac{\delta(t-\tau)}{1-\int_{0}^{\tau -x}f_{S}(t')\,dt'}\int_{t-x}^{\infty}f_{S}(\tau)\,d\tau, \label{SX sharp sum}
\end{align}
where $\theta$ is the Heaviside step function. The expectations of the above random variables are then given by
\begin{align}
&\mathbf{E}[S_{+}(x)] = \int_0^{\infty} t f_{S_{+}}(t) dt = \frac{\int_{0}^{\tau - x} t f_{S}(t) dt}{\int_{0}^{\tau - x}f_{S}(t')\,dt'}, \label{SX sharp sum mean} \\
&\mathbf{E}[R_{+}(x)]=\int_0^{\infty} t f_{R_{+}}(t) dt = \frac{\tau \left( 1 - \int_{0}^{\tau - x}f_{S}(t')\,dt' \right)}{1 - \int_{0}^{\tau - x}f_{S}(t')\,dt'} = \tau. \label{RX sharp sum mean}
\end{align}
Substituting the expectation for $\mathbf{E}[S_{+}(x)]$ and $\mathbf{E}[R_{+}(x)]$ given above into Eq. (\ref{General Mean +}) yields the following formula for the mean service time under sharp resetting
\begin{align}
\mathbf{E}[B_{\tau}] = \, \mathbf{E}[X] + \int_{0}^{\infty} \frac{\int_{0}^{\tau - x} t f_{S}(t) dt}{\int_{0}^{\tau - x}f_{S}(t')\,dt'} f_X(x) \, dx  +\tau \left (\int_{0}^{\infty} \frac{f_X(x)}{\int_{0}^{\tau - x}f_{S}(t')\,dt'} \, dx -1 \right ). \label{Sharp Mean +}
\end{align}

We emphasize that the results derived so far in Secs. \ref{Sec multiplicative} and \ref{Sec additive} are agnostic to the details of the arrival process and make no assumptions about its stochastic behavior. In the following section we consider the celebrated M/G/1 queue for which we derive explicit expressions for the probability generating function of the queue length and for the Laplace transform of the waiting times, as well as for the mean queue size.

\section{Queue length distribution of $S\&X$-M/G/1 queues with service resetting}
\label{Sec distribution}

So far, we were concerned only with the mean service time under resetting. We now turn our attention to the distribution of the queue length. To this end, we consider an M/G/1 queue with Poisson arrivals with rate $\lambda$ and overall service time $U$. In an M/G/1 queue, the probability mass function of the number of jobs in the system, $L$, is denoted as $P_L(n)=\text{Pr}(L=n)$. The probability generating function (PGF) 
\begin{align} 
 G_{L}(z)  = \sum_{n=0}^{\infty} P_{L}(n)z^{n}~, \label{PGF def}
\end{align}
of the queue length, is given by the Pollaczek–Khinchine formula \parencite{Kleinrock-book}
\begin{align}
G_{L} (z)=\frac{(1-\rho)(1-z)\tilde{U}(\lambda (1-z))}{\tilde{U}(\lambda (1-z))-z}~, \label{PK transform}
\end{align}
where $\tilde{U}(\cdot)$ is the Laplace transform of the service time, and $\rho=\lambda \mathbf{E}[U] < 1$ is the utilization of the queue. The probability that there are $n$ jobs in the system can then be computed from Eq. (\ref{PK transform}) via
\begin{equation}
P_{L}(n)=\text{Pr}(L=n)=\frac{G_{L}^{(n)}(0)}{n!}~, \label{PGF inverse}
\end{equation}
where $G_{L}^{(n)}(0)$ stands for the n-th derivative of $G_{L}(z)$ evaluated at $z=0$.
Similarly, one can also compute the Laplace-Stieltjes transform, $\tilde{W}(s)$, of the total time a job spends in the system.  This is given by the corresponding Pollaczek-Khinchin transform equation which reads
\begin{equation}
\tilde{W}(s)=\frac{\tilde{U}(s)(1-\rho)s}{s-\lambda(1-\tilde{U}(s))}~. \label{Sojourn time}
\end{equation}
The mean queue length is given by 
\begin{equation}
\mathbf{E}[L]=\frac{\rho}{1-\rho}+\frac{\rho^2}{2(1-\rho)}\left( CV_{U}^{2}-1  \right)~,
\label{PK-1}
\end{equation}
where $CV_{U}^{2}$ is the squared coefficient of variation of the overall service time $U$, $CV_{U}^2=\sigma^2(U)/\mathbf{E}^2[U]$.

To obtain corresponding formulas in the case of service resetting, we substitute the Laplace transforms of the service time with resetting into the Pollaczek–Khinchine formula. For example, substituting   $\tilde{V}_r(s)$ from Eq. (\ref{Main-Multiplication-Exp}) and $\tilde{B}_r(s)$ from Eq. (\ref{Main-Additive-Exp}) into the Pollaczek–Khinchine formula yields the PGF of the queue length under Poissonian resetting, in the multiplicative and additive models, respectively.

\begin{figure*}[t!]
\begin{centering}
\includegraphics[width=1\textwidth]{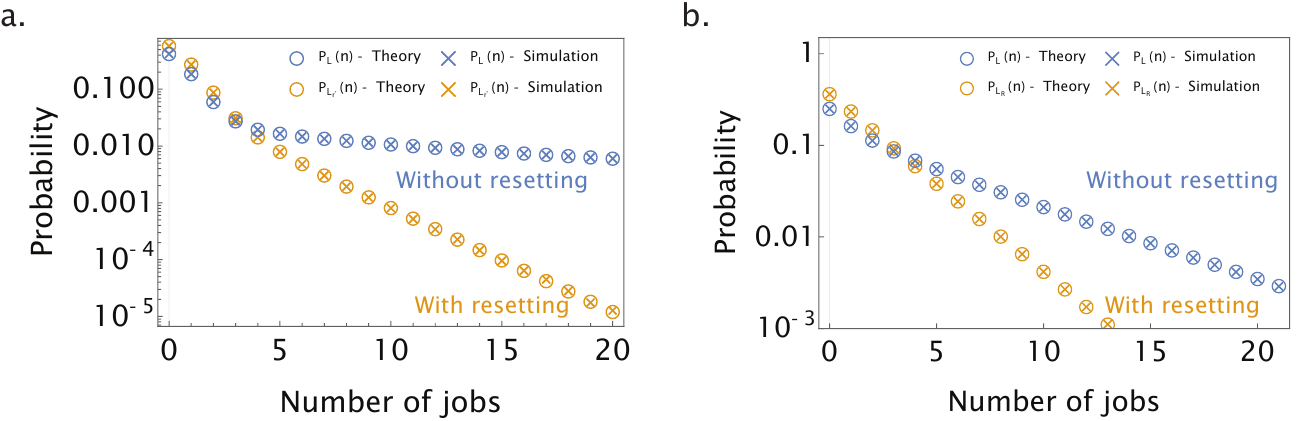}
\end{centering}
\caption{The probability mass function of the number of jobs in the $S\&X$-M/G/1 queueing system (both in queue and in service) without service resetting (blue) and with optimal Poissonian service resetting (orange). The probabilities $P_{L_{r^*}}(n)$ and $P_{L}(n)$ to find $n$ jobs in a system with and without service resetting, respectively, were computed by use of Eq. (\ref{PGF inverse}). Here, $r^*$ was taken as the optimal resetting rate which minimizes the overall mean service time. The analytical results (circles) were further corroborated with numerical simulations (X marks). Panel (a): The additive model. Here, the arrival rate is set to $\lambda = 1/2$ and the service time $B$ is taken as the sum of: (i) $S$, a Gamma random variable with parameters $\alpha = 0.01$ and $\theta = 50$, and (ii) $X$, a deterministic random variable $\text{Pr}(X=2/3)=1$. The optimal resetting rate is then found to be $r^* \simeq 0.2424$ (see Fig. \ref{Fig2}a). Panel (b): The multiplicative model. Here, the jobs arrival rate is set to $\lambda = 3/4$, and the service time $V$ is taken to be a product of: (i) $S$, an Inverse Gaussian random variable with parameters $\alpha_S = 3/2$ and $\beta_S = 3/4$, and (ii) $X$, a deterministic random variable $\text{Pr}(X=2/3)=1$. The optimal resetting rate is then found to be $r^* \simeq 1.372$ (see Fig. \ref{Fig4}a). In both the additive and multiplicative models, it is evident that the introduction of optimal service resetting gives rise to lighter-tailed queue length distributions and consequently to shorter queues.}
\label{Fig6}
\end{figure*}

In Fig. \ref{Fig6}, we follow the aforementioned steps to derive the queue length distribution under resetting for both the additive and multiplicative models. In Fig. \ref{Fig6}a, we consider the additive model and compare the queue length distribution under optimal Poissonian resetting to the queue length distribution without resetting. We take the distribution of the server slowdown to be Gamma with $\mathbf{E}[S]=1/2$ and $\sigma^2(S)=25$, and fix the inherent job size to $x=2/3$ (i.e., identical to the $\alpha=0.01$ case from Fig. \ref{Fig2}a). The optimal resetting rate in this scenario is $r^* \simeq 0.2424$. The effect of resetting in this scenario is evident --- under resetting, there is a higher probability of having shorter queues compared to the no  resetting case. Furthermore, the tail of the distribution for a queue with resetting is significantly lighter. The effect is even more pronounced when comparing the mean queue lengths: the mean queue length under optimal resetting is $\mathbf{E}[L_{r^*}] \simeq 0.739$, which is more than ten fold smaller than its value $\mathbf{E}[L] \simeq 8.49$ without resetting. 

%%%%%%%%%%%%%%%%%%%%%%%%%%%%%%%%%%%%%%%%%%

%%%%%%%%%%%%%%%%%%%%%%%%%%%%%%%%%%%%%%%%%%

In Fig. \ref{Fig6}b, we consider the multiplicative model and compare the queue length distribution under optimal Poissonian resetting to the queue length distribution without resetting. Here, we take the distribution of the server slowdown to be Inverse Gaussian with $\mathbf{E}[S]=3/2$ and $\sigma^2(S)=9/2$, and fix the inherent job size as $x=2/3$ (i.e., identical to the $\beta_S=0.75$ case from Fig. \ref{Fig4}a). The optimal resetting rate in this scenario is $r^* \simeq 1.372$. The effect of resetting under the multiplicative model is qualitatively similar to the additive case. Under resetting, there is a higher probability of having shorter queues compared to the no resetting case. Once again, the tail of the queue length distribution with resetting is lighter. We find that optimal resetting reduces the mean queue length from $\mathbf{E}[L] \simeq 4.12$ to $\mathbf{E}[L_{r^*}] \simeq 1.76$, i.e., more than two-fold.

Summarizing, in this section we see the effects of optimal Poissonian resetting on the queue length distribution of an M/G/1 queue in both the additive and multiplicative models. In both cases, we demonstrated that resetting has a significant impact on the queue length distribution: generating lighter tails, smaller means, and increasing the probability of seeing shorter queues.

\section{Conclusions}
\label{Conclusions}
In this paper, we further developed the theory of queueing systems with service resetting, extending to $S\&X$ queues \parencite{S&Xmodel} the framework presented in \parencite{Resettingqueues}. We considered scenarios where a job's service time is either a product, or a sum, of two components: one intrinsic to the server, the other extrinsic. We derived the distribution and mean of the total service time of a job in cases where times between service resetting events form a general renewal process. The prevalent constant-rate (Poissonian) resetting policy was analyzed as an example, and explicit conditions under which resetting reduces the mean service time and improves queue performance were derived. In the multiplicative case, the condition is given by Eq. (\ref{CV condition X}), which involves only the mean and variance of the intrinsic component of the service time. In the additive case, the condition is given by Eq. (\ref{CV condition +}) which involves both the means and variances of the intrinsic and extrinsic components of the service time. In both cases, the introduction of resetting is beneficial when the variance of the intrinsic component of the service time is large enough. Aside from Poissonian resetting, we also analyzed sharp (deterministic) resetting, establishing results for cases where times between resetting events are fixed.

In our analysis, we have made the common assumption that the arrival process to the queue is independent of the service process. Thus, the results obtained in Secs. \ref{Sec multiplicative} and \ref{Sec additive} hold regardless of the arrival process. For illustration, we considered the canonical M/G/1 queue for which we have given explicit results for the distribution and mean of the queue length in $S\&X$ queueing systems with service resetting. Yet, we  emphasize that results and conclusions coming from our analysis carry over to other systems such as $S\&X$-G/G/1 queues. Indeed, the conditions specified in Eq. (\ref{CV condition X}) and Eq. (\ref{CV condition +}) guarantee that service resetting will reduce the mean service time irrespective of the arrival process.

In summary, this paper established the pivotal role of service resetting as a novel strategy for combating the harmful effects of service time fluctuations in $S\&X$ queues. Our comprehensive analysis, spanning both additive and multiplicative models, highlights the effectiveness of service resetting in significantly reducing mean service times and, consequently, shrinking queue lengths. The developed framework not only offers theoretical insights but also equips practitioners with a practical toolkit. By providing explicit conditions for the implementation of beneficial resetting, our results enable precise decision-making, ensuring that service resetting becomes an invaluable asset in queue management scenarios. Notably, our findings transcend the confines of specific arrival processes, underscoring the widespread applicability of service resetting.

\section{Appendix}

\subsection{Derivation of Eq. (\ref{mean service time expansion X})}
\label{Appendix B}
\noindent Equation (\ref{mean service time expansion X}) can be derived by expanding Eq. (\ref{Mean-Multiplication-Exp}) in a Taylor series to first order in $\delta r$, where $\delta r$ is an infinitesimally small resetting rate. We start by expanding both the numerator and denominator of Eq. (\ref{Mean-Multiplication-Exp}). Utilizing the moment expansion of the Laplace transform, one can expand the numerator of the right hand side of Eq. (\ref{Mean-Multiplication-Exp}) to second order in $\delta r$. This yields 
\begin{align}
\tilde{S}(\delta r x) =  ~ 1 -\delta r x \mathbf{E}[S] + \frac{1}{2}(\delta r x)^2 \mathbf{E}[S^2]+ O(\delta r^3). \label{Srx LT}
\end{align}
Following the above steps, we obtain the Taylor series of $r\tilde{S}(rx)$ to second order in $\delta r$
\begin{align}
\delta r\tilde{S}(\delta rx) =  ~ \delta r -\delta r^2 x \mathbf{E}[S] + O(\delta r^3). \label{rSrx LT}
\end{align}
Substituting Eqs. (\ref{Srx LT}) and (\ref{rSrx LT}) into Eq. (\ref{Mean-Multiplication-Exp}),  and dividing both numerator and denominator by $\delta r$ yields
\begin{equation}
\mathbf{E}[V_{\delta r}] = ~ \int_{0}^{\infty}\frac{x \mathbf{E}[S] - \frac{1}{2}\delta r x^2 \mathbf{E}[S^2]+ O(\delta r^2)}{1 -\delta r x \mathbf{E}[S] + O(\delta r^2)}f_{X}(x)\,dx. \label{Mean-Multiplication-Exp-der1}
\end{equation}
Using the approximation $\frac{1}{1-\epsilon} \simeq 1 + \epsilon $ where $\epsilon \rightarrow 0$, Eq. (\ref{Mean-Multiplication-Exp-der1}) reduces to 
\begin{align}
\mathbf{E}[V_{\delta r}] = \int_{0}^{\infty}&\Big[ x \mathbf{E}[S] - \frac{1}{2}\delta r x^2 \mathbf{E}[S^2]+ O(\delta r^2) \Big] \Big[ 1 + \delta r x \mathbf{E}[S] + O(\delta r^2) \Big] f_{X}(x)\,dx. \label{Mean-Multiplication-Exp-der2}  
\end{align}
Expanding the product given in the above equation and retaining terms up to first order in $\delta r$ yields
\begin{align}
\mathbf{E}[V_{\delta r}] = \int_{0}^{\infty}&\Big( x \mathbf{E}[S] + \delta r x^2 \left( \mathbf{E}[S]^2 - \frac{1}{2}\mathbf{E}[S^2] \right) + O(\delta r^2) \Big)f_{X}(x)\,dx. \label{Mean-Multiplication-Exp-der3}  
\end{align}
Integrating with respect to $x$, Eq. (\ref{mean service time expansion X}) readily follows.

\subsection{Derivation of Eq. (\ref{mean service time expansion +})}
\label{Appendix A}
We derive Eq. (\ref{mean service time expansion +}) by expanding both numerator and denominator of Eq. (\ref{Mean-Additive-Exp}) in a Taylor series. Utilizing the moment expansion of the Laplace transform, one can expand the numerator of Eq. (\ref{Mean-Additive-Exp}) to second order in $\delta r$
\begin{align}
\tilde{X}(-\delta r)-\tilde{S}(\delta r) & = ~ \Big(1 + \delta r\mathbf{E}[X] + \frac{1}{2}\delta r^2 \mathbf{E}[X^2] \Big) - \Big(1 - \delta r\mathbf{E}[S] + \frac{1}{2}\delta r^2 \mathbf{E}[S^2] \Big) +O(\delta r^3)  \nonumber \\
&=\delta r \Big(\mathbf{E}[X] + \mathbf{E}[S] \Big) + \frac{1}{2} \delta r^2 \Big(\mathbf{E}[X^2] - \mathbf{E}[S^2] \Big) + O(\delta r^3). \label{numerator + expansion}
\end{align}
Similarly, we expand the denominator of the right hand side of Eq. (\ref{Mean-Additive-Exp}) to second order in $\delta r$
\begin{align}
\delta r\tilde{S}(\delta r) = \delta r - \delta r^2 \mathbf{E}[S] + O(\delta r^3). \label{denominator + expansion}
\end{align}
Substituting Eqs. (\ref{numerator + expansion}) and (\ref{denominator + expansion}) into Eq. (\ref{Mean-Additive-Exp}), and dividing both numerator and denominator by $\delta r$, yields
\begin{equation}
\mathbf{E}[B_{\delta r}] = \frac{\mathbf{E}[X] + \mathbf{E}[S] + \frac{1}{2} \delta r \Big(\mathbf{E}[X^2] - \mathbf{E}[S^2] \Big) + O(\delta r^2)}{1 - \delta r \mathbf{E}[S] + O(\delta r^2)}. \label{Mean-Additive-Exp expansion 1}  
\end{equation}
Using the approximation $\frac{1}{1-\epsilon} \simeq 1 + \epsilon $ where $\epsilon \rightarrow 0$, Eq. (\ref{Mean-Additive-Exp expansion 1}) reduces to 
\begin{align}
\mathbf{E}[B_{\delta r}] = &\Big[ \mathbf{E}[X] + \mathbf{E}[S] + \frac{1}{2} \delta r \Big (\mathbf{E}[X^2] - \mathbf{E}[S^2] \Big) + O(\delta r^2) \Big ]  \Big[ 1 + \delta r \mathbf{E}[S] + O(\delta r^2) \Big]. \label{Mean-Additive-Exp expansion 2}  
\end{align}
Expanding the product given in the above equation and retaining terms up to first order in $\delta r$ yields
\begin{align}
\mathbf{E}[B_{\delta r}] = & ~\mathbf{E}[X] +\mathbf{E}[S] +\delta r \Big( \mathbf{E}[X] \mathbf{E}[S] + \mathbf{E}[S]^{2} + \frac{1}{2}\left(\mathbf{E}[X^2] -\mathbf{E}[S^2] \right)\Big) +O(\delta r^2), \label{Mean-Additive-Exp expansion 3}
\end{align}
which is Eq. (\ref{mean service time expansion +}) in the main text. 

\section{Acknowledgments}
\noindent This project has received funding from the European Research Council (ERC) under the European Union’s Horizon 2020 research and innovation program (grant agreement No. 947731). U.Y acknowledges support from the Israel Science Foundation (grant No. 1968/23).

\printbibliography

\end{document}